\documentclass{article}

\newcommand{\bigO}{\mathcal{O}}
\newcommand{\RNum}[1]{\uppercase\expandafter{\romannumeral #1\relax}}  %%roman numerals
\newcommand{\abs}[1]{\lvert#1\rvert}
\newcommand{\nameA}{Simple Sweep}
\newcommand{\nameB}{Window-wise Sweep}
\newcommand{\nameC}{Corrective Sweep}

\usepackage{amssymb}
\usepackage{amsmath}
\usepackage[numbers]{natbib}
 \usepackage{hyperref} % move hyperref package to the beginning to
\usepackage[utf8]{inputenc}
\usepackage[T1]{fontenc}
\usepackage[british]{babel}
\usepackage{mathptmx}       % selects Times Roman as basic font
\usepackage{helvet}         % selects Helvetica as sans-serif font
\usepackage{courier}        % selects Courier as typewriter font
\usepackage{type1cm}        % activate if the above 3 fonts are
\usepackage{csquotes}
\usepackage{etex}
\usepackage{caption}
\usepackage[textsize=tiny]{todonotes}
\usepackage{url}
\usepackage{comment}
\usepackage{algorithm}
\usepackage{algpseudocode}
\usepackage{tikz}
\usetikzlibrary{calc}
\usetikzlibrary{decorations.pathmorphing,patterns}
\usetikzlibrary{calc,patterns,decorations.markings}
\usetikzlibrary{decorations.pathreplacing}
\usepackage{tikz}
\usepackage{paralist}
\usetikzlibrary{calc}
\usetikzlibrary{decorations.pathmorphing,patterns}
\usetikzlibrary{calc,patterns,decorations.markings}
\usetikzlibrary{decorations.pathreplacing}
\usepackage{enumerate}
\usepackage{siunitx} %format numbers in table
\usepackage{array}
\usepackage{multirow}
\usepackage{subfig}
\usepackage{geometry}
\usepackage[margin=1cm]{caption}
\usepackage{paralist}
\usepackage{fancyhdr}

\providecommand{\keywords}[1]{\textbf{\textit{Key words.}} #1}

\providecommand{\keywords}[1]{\textit{Keywords:} #1}

\setcitestyle{square, numbers}

\begin{document}

\title{Sweep Algorithms for the Capacitated  Vehicle Routing Problem with
 Structured Time Windows} % $^{\dag}$}
% your contribution title if the original one is too long

\author{
C. \ Hertrich \thanks{Department of Mathematics, Technische Universit\"at Kaiserslautern, Germany,
\href{mailto:hertrich@math.tu-berlin.de}{hertrich@math.tu-berlin.de}
\newline
Christoph Hertrich is supported by the Karl Popper Kolleg
“Modeling-Simulation-Optimization”
funded by the Alpen-Adria-Universit\"at Klagenfurt and the
Carinthian Economic
Promotion Fund (KWF).
}
\ ,
P.\ Hungerländer\thanks{Department of
Mathematics, Alpen-Adria Universität Klagenfurt, Austria,
  \href{mailto:philipp.hungerlaender@aau.at}{philipp.hungerlaender@aau.at}}
 \ , and C.\
  Truden\thanks{Department of
  Mathematics, Alpen-Adria Universität Klagenfurt, Austria,
  \href{mailto:christian.truden@aau.at}{christian.truden@aau.at}}
  }

\maketitle

\abstract{%
The capacitated Vehicle Routing Problem with structured Time Windows (cVRPsTW)
 is concerned with finding optimal tours for vehicles with given capacity
 constraints to deliver goods to customers within assigned time windows.
 In our problem variant these time windows have a special structure,
 namely they are non-overlapping and each time window holds several customers.
 This is a reasonable assumption for Attended Home Delivery services.
Sweep algorithms are known as simple, yet effective
heuristics for the classical capacitated Vehicle Routing Problem.
We propose variants of the sweep algorithm that are not only able to deal
with time windows, but also
exploit the additional structure of the time windows in a cVRPsTW.
Afterwards we suggest local
improvement heuristics to decrease our objective function even further.
A carefully constructed benchmark set that resembles real-world data is used
to prove the
efficacy of our algorithms in a computational study.
}
\\

\maketitle

\section{Introduction}

\text{Attended Home Delivery} (AHD)
services, e.g., online grocery shopping services,
have encountered a significant growth in popularity in recent years.
In a typical setup,
customers choose time windows during which they want to receive their
goods from a set of available time windows that is provided by the supplying company.
This set is called \textit{structured} if the number of customers is
significantly larger  than the number of time windows
and  all windows are pair-wise non-overlapping.
Once all orders have been placed, the supplier aims to minimize the fulfillment
costs. This involves solving a
\textit{capacitated Vehicle Routing Problem with structured Time Windows} (cVRPsTW),
which is a special case of the \emph{capacitated Vehicle Routing
Problem with Time Windows} (cVRPTW).
Heuristics are the method of choice in practice
 to produce high-quality solutions in
reasonable time since the problem instances occurring are typically rather large.

So-called \textit{cluster-first, route-second} methods have
been proven to be effective for the classical
\emph{capacitated Vehicle Routing Problem} (cVRP).
These methods first partition the customers into subsets
that are small enough such that they can all be visited by one vehicle.
In a second step, they compute a route for each vehicle.
The most prominent example of this group of algorithms is the
\textit{sweep algorithm}
 \cite{Gillet1974}.
 It clusters the customers by dividing
the plane into radial sectors originating from the depot's location.
In this work we first generalize the sweep algorithm
to the cVRPTW.
Secondly, we propose two variants that exploit
the additional structure of the time windows in a cVRPsTW.
 Due to the imposed structure the
total number of time windows is quite small, which allows us to use
window-dependent angles.
 After obtaining a first feasible solution,
improvement heuristics are applied that try to decrease the objective
function by slightly altering the angles obtained by our sweep algorithms.

  Finally, we conduct an extensive computational study
  using a large variety of carefully constructed benchmark instances
 which show that
our approach is capable of finding good initial solutions for instances
containing up to 2000 customers within a few seconds.
Further, we demonstrate that the
proposed improvement heuristics allow us
to significantly  improve the solution
quality within a few minutes.
We notice that the performance of the different variants is dependent
on the characteristics of the considered instance, e.g., whether vehicle
capacities or time windows are the stronger restriction.

For an overview of exact  resp. heuristic  methods for the cVRP(TW) we refer to
\cite{BALDACCI20121,VRPTW_ReviewI}.
Solomon \cite{Solomon:1987} proposes a sweep heuristic that
takes time window constraints into account.
In contrast to our approach,  it
does not consider the time windows of the customers when partitioning them.
Instead, the time windows are only respected when computing the routes for each vehicle.
While Solomon's approach utilizes an insertion heuristic to obtain the routes,
 we apply a Mixed-Integer Linear Program (MILP)
 to decide the feasibility of an
  assignment of customers to a vehicle,
as well as to
obtain the optimal solution of the
occurring routing subproblems.

\section{Formal Problem Definition}\label{description}
A cVRPsTW instance consists of a set of time windows
$\mathcal{W}=\{w_1,\dots,w_q\}$, where each $w_i\in\mathcal{W}$
is defined through its start time $s_{w_i}$ and its end time $e_{w_i}$ with
$s_{w_i}<e_{w_i}$,
a set of customers
$\mathcal{C}$, $\abs{\mathcal{C}}=n$, a time window assignment function
$w\colon \mathcal{C}\to \mathcal{W}$,
a depot $d$ from which all vehicles depart from and return to,
$\overline{\mathcal{C}} := \mathcal{C}\cup\{d\}$,
a travel time function
$t\colon \overline{\mathcal{C}}\times\overline{\mathcal{C}}\to \mathbb{R}_{\geq0}$,
a service time function $s\colon \mathcal{C}\to \mathbb{R}_{>0}$,
a common vehicle capacity $C\in\mathbb{R}_{>0}$, and an order weight function
$c\colon \mathcal{C}\to \left]0,C\right]$. We require the windows to be
structured, i.e.,
$n\gg q$ and for $1\leq i < j \leq q$ it holds $s_{w_j}\geq e_{w_i}$.
Moreover, each element in $\overline{\mathcal{C}}$ has
coordinates
 in the two-dimensional plane.
We assume that the travel times are correlated to the geographical distances, but
not purely determined by them.

A \textit{tour} consists of a set
$\mathcal{A}=\{a_1,a_2,\dots, a_k\}$ of customers with corresponding arrival
times $\alpha_{a_1},\dots,\alpha_{a_k}$
during which the vehicles are scheduled to arrive. A tour is called
\textit{capacity-feasible}, if $\sum_{i=1}^{k} c(a_i) \leq C$.
Moreover, we call it \textit{time-feasible}, if every customer
is served within its assigned time window,
i.e., $s_{w(a_i)}\leq\alpha_{a_i}\leq e_{w(a_i)}$,
$i=1,\ldots,k$, and if there is sufficient time to respect the required
service and travel times, i.e.,
$\alpha_{a_{i+1}}-\alpha_{a_i}\geq s(a_i)+t(a_i,a_{i+1})$,
$i=1,\ldots,k-1$.
A \textit{schedule} $\mathcal{S} = \{\mathcal{A},\mathcal{B},\dots\}$ is
a set of tours where each customer occurs in exactly one tour.
It is called \textit{feasible}, if all tours are capacity- and time-feasible.

We consider three objectives:
The first one is the
\textit{number  of vehicles used}, i.e.\
$\lambda_1(\mathcal{S}):=\abs{\mathcal{S}}$.
Secondly, the \emph{schedule duration} $\lambda_2(\mathcal{S})$
is defined as the sum of all \emph{tour durations}, i.e.\
$\lambda_2(\mathcal{A}) :=t(d,a_1) + \alpha_{a_k}-\alpha_{a_1} +
s(a_k) + t(a_k,d)$.
Thirdly, the \emph{schedule travel time} $\lambda_3(\mathcal{S})$
is defined as the sum of all \emph{tour travel times}, i.e.\
$\lambda_3(\mathcal{A}) :=t(d,a_1) + \sum_{i=1}^{k-1} t(a_i,a_{i+1}) + t(a_k,d)$.
Similar to
Solomon \cite{Solomon:1987},
 we aim to minimize these three objectives
with respect to the lexicographical order $(\lambda_1, \lambda_2, \lambda_3)$,
since providing a vehicle is usually the most expensive cost component,
followed by the drivers' salaries, and the costs for fuel.

\section{Sweep Algorithms for Structured Time Windows}\label{sweep}

In this section we describe several  variants of the sweep algorithm.
First, we introduce some more notation and definitions.

\textbf{Tree-Feasibility:}
For a
tour $\mathcal{A}$
let  $\mathcal{A}_{w}\subseteq\mathcal{A}$
be the set of customers assigned to $w \in \mathcal{W}$.
We consider the complete directed graph
with vertex set $\mathcal{A}_{w}$ and edge weights $t(a_i, a_j) + s(a_i)$
assigned to each edge $(a_i,a_j)$.
The existence of a spanning arborescence through $\mathcal{A}_w$ with length of at
most $e_w- s_w$ for all time windows $w$ forms a necessary condition
for the time-feasibility of a tour
% through
 $\mathcal{A}$.
 In this case we call $\mathcal{A}$
\textit{tree-feasible}.
Existence of such an arborescence can be checked in $\bigO(|\mathcal{A}_w|^2)$ time \cite{Tarjan1977}.
Hence,
time-infeasibility of a potential tour can often
be detected without solving a time-consuming MILP.
Moreover, we apply  the concept of tree-feasibility in the \nameC \
algorithm below.

\textbf{Angles of Customers:}
Any sweep algorithm is based on the polar coordinate representation of the
customers $\mathcal{C}$,
where the depot $d$ forms the origin of the coordinate system and
$\theta(a)\in[0,2\pi[$ denotes the angle component of a customer $a \in \mathcal{C}$.
The direction of the zero angle $\theta_0$ is a choice parameter of the algorithm as it
impacts the result.
We choose the zero angle such that it separates
the two consecutive customers with the largest angle gap, i.e.,
$\max_{a,b\in\mathcal{C}}\theta(a)-\theta(b)$ is minimized.
Moreover, we run all algorithms in clockwise and counterclockwise direction
and in each case select the variant that gives the  better result.

In the following, we denote
$\mathcal{C}(\theta,\theta')=\{a\in\mathcal{C}\mid\theta(a)\in[\theta,\theta'[\}$,
 $\mathcal{C}_w=\{a\in\mathcal{C}\mid w(a)=w\}$,
 and $\mathcal{C}_w(\theta,\theta')=\mathcal{C}(\theta,\theta')\cap\mathcal{C}_w$.

\subsection{Sweep Strategies}
We propose the following general strategy to obtain high-quality
solutions for a given
cVRPTW or cVRPsTW instance. It consists of three steps:

\begin{compactenum}
  \setlength{\itemsep}{0pt}
  \item \label{step:clustering} Use a variant of the  sweep algorithm
    to determine a feasible clustering of $\mathcal C$.
  \item \label{step:localImpr}Apply local improvement heuristics to enhance the quality of the clustering.
  \item \label{step:routing} Compute the optimal route for each cluster.
\end{compactenum}

In all three steps, the \textit{Traveling Salesperson Problem with
  Time Windows} (TSPTW) respectively the
\textit{Traveling Salesperson Problem with structured Time Windows} (TSPsTW)
occurs as a
subproblem to check the  time-feasibility or to obtain the optimal solution of a
single tour.
We apply two MILP formulations that have been proposed in previous
work \cite{Hungerlaender201801}, a general one for the TSPTW,
and a more efficient one that is tailored to the TSPsTW.
Following the lexicographical order, it first minimizes
 $\lambda_2$,  and then $\lambda_3$ while keeping $\lambda_2$ fixed.
Next let us relate our notation to the well-known
sweep algorithm \cite{Gillet1974} for the cVRP.

 \textbf{Traditional Sweep:}
The clustering method proposed in \cite{Gillet1974}
relates to Step \ref{step:clustering}.
It finds angles $\theta_i$ such that the $i$-th cluster is given
by $\mathcal{C}(\theta_{i-1},\theta_{i})$ as follows:
Set $\theta_0=0$.
For $i=1,2,\ldots$
make $\theta_i$ as large as possible such that
$\mathcal{C}(\theta_{i-1},\theta_{i})$ forms a capacity-feasible cluster.
Here, and in all following algorithms, the range for $i$ is chosen
such that all customers are scheduled.

Next, we propose a natural generalization of the Traditional Sweep algorithm that
works for instances having structured (cVRPsTW), as well as for instances having
 arbitrary time windows (cVRPTW).

\textbf{\nameA:}
Choose the angle   $\theta_i$ as large as possible while ensuring that the resulting cluster
$\mathcal{C}(\theta_{i-1},\theta_{i})$ is still small enough such that
a time- and capacity-feasible tour that visits all contained customers can be found.
We check the time-feasibility using the TSPTW- resp. TSPsTW-MILP.

Now we present two variants of the sweep algorithm that exploit
the additional structure of the time windows of cVRPsTW instances.

\textbf{\nameB:}
In case of structured time windows
 there are quite few time windows in
comparison to customers, i.e., $q \ll n$.
Hence, we can define window-dependent angles
$\theta_{i}^{j}, ~i \geq 1, ~j=1,\ldots,q,$ such that
the $i$-th cluster is given by
$\bigcup_{j=1}^q \mathcal{C}_{w_j}(\theta_{i-1}^{j},\theta_i^{j})$.
We propose to add the customers to the clusters window by window.
The resulting algorithm is described as follows:
 While there are
unclustered customers in $\mathcal{C}_{w_j}$, make $\theta_{i}^{j}$ as
large as possible such that $\mathcal{C}_{w_j}(\theta_{i-1}^{j},\theta_i^{j})$
can be added to the $i$-th cluster while ensuring that the cluster
stays  time- and capacity-feasible. If necessary, increase the number of clusters.

\textbf{\nameC:}
If vehicle capacity is a stronger restriction than time windows, then \nameB\ creates
large sectors for the first few time windows and runs out of capacity later on,
causing an increased need of vehicles.
Therefore, we propose another variant that prevents this from happening.
We first initialize angels
 $\theta_{i}^{0}$ of maximal size such that all sets
$\mathcal{C}(\theta_{i-1}^{0},\theta_{i}^{0})$
are capacity- and tree-feasible.
Then we start with empty clusters. For $j=1,\ldots,q$ we set
$\theta_{i}^{j}=\theta_{i}^{j-1}$ and add the customers of
$\mathcal{C}_{w_j}$ accordingly.
Since this may result in some infeasible
clusters, the angles $\theta_{i}^j$ have to be adjusted.
Starting with $i=1$, we check whether
the $i$-th cluster is still feasible.
 If not, we try to reduce the size of cluster $i$ by increasing the angles
$\theta_{i-1}^{j},\dots,\theta_{1}^{j}$, such that
cluster $i$ becomes feasible  while the clusters $1,\dots,i-1$ remain feasible.
If this procedure does not succeed,
we decrease $\theta_{i}^{j}$ until the $i$-th cluster becomes feasible.
Then, we increment $i$ and repeat this procedure.
If necessary, we increase the number of clusters.

\textbf{Local Improvement Heuristic:}
After obtaining an initial clustering using one of the
heuristics described above, we aim to improve the clustering during Step
\ref{step:localImpr}.
In case that the clustering was obtained using the \nameA\ algorithm,
we set $\theta_{i}^{j}=\theta_{i}$.
For each angle $\theta_{i}^{j}$ we try to improve the objective function
iteratively by decreasing or increasing the angle. In each step the angle is
slightly altered such that one customer moves to another cluster.
A change is accepted if the lexicographical objective is improved.
This procedure is repeated until a local minimum is reached.

\textbf{Routing:}
In Step \ref{step:routing} a tour for each cluster is obtained
by solving the TSPTW- or TSPsTW-MILP with lexicographical objective $(\lambda_2,\lambda_3)$ to optimality.

\section{Computational Experiments}\label{experiments}

As to our best knowledge, none of the available cVRPTW benchmark instances
comply with the considered AHD use-case,
we created
a new benchmark set that resembles urban settlement structures
in order to provide meaningful computational experiments.
We placed customers on a 20\,km $\times$ 20\,km square grid
(roughly the size of Vienna).
In order to achieve varying customer densities,
only 20\,\% of the customer locations are sampled from a two-dimensional
uniform distribution.
The remaining 80\,\% of the customer locations have been sampled from
clusters whose centers and shapes have been randomly sampled.
Travel times are calculated proportional to the Euclidean distances rounded to integer seconds.
As proposed by Pan et al. \cite{pan2017}
 we assume a travel speed of $20$\,km/h.
Each customer is randomly assigned to one
time window out of a set of ten consecutive time windows, where each is one hour long.
We assume the service time at each customer to be 5 minutes.
The order weights are sampled
from a truncated normal distribution centered around 5 units.
We consider
configurations with $n=250,500,1000,2000$  and   $C=200,400$.
Choosing  $C=200$ results in the  vehicle capacity being the bottleneck,
while choosing  $C=400$ causes that the time window constraints are the
 most limiting factor.
For each configuration we constructed $100$
instances and report the average results.
All instances can be downloaded from \url{http://tinyurl.com/vrpstw}.

All experiments were performed on an Ubuntu 14.04 machine equipped
with an Intel Xeon E5-2630V3 @ 2.4\,GHz 8 core processor and 132\,GB RAM.
We use Gurobi 8.0 in single thread mode to solve the MILPs.
To demonstrate the effectiveness of
the improvement heuristics, we compare the case where all three steps are performed
against
the case where Step \ref{step:localImpr} is omitted.

\setlength{\tabcolsep}{4pt}
\renewcommand{\arraystretch}{0.84}
\begin{table}[]
\centering
\begin{tabular}{|l|rrrr|rrrr|}
	\hline
	Vehicle Capacity & \multicolumn{4}{c|}{200} & \multicolumn{4}{c|}{400} \\
	Runtime/Objectives  & $t$ & $\lambda_1$  & $\lambda_2$ & $\lambda_3$ & $t$&  $\lambda_1$ & $\lambda_2$ & $\lambda_3$ \\
	\hline
	\textbf{Without Improvement}&&&&&&&&\\
	\nameA 			& 7.2 & 54.0 & 507.7 & 332.7 & 41.0 & 42.1 & 406.5 & 315.3 \\
	\nameB	 		& 96.9 & 64.4 & 381.2 & 371.5 & 306.7 & 35.0 & 340.2 & 328.1 \\
	\nameC		 	& 7.6 & 54.0 & 507.2 & 332.5 & 34.6 & 41.5 & 402.1 & 314.8 \\
	\hline
	\textbf{With Improvement}&&&&&&&&\\
	\nameA 			& 241.0 & 54.0 & 480.3 & 318.1 & 1324.7 & 41.9 & 368.5 & 300.8 \\
	\nameB	 		& 141.9 & 64.4 & 378.5 & 369.3 & 679.1 & 35.0 & 336.2 & 324.8 \\
	\nameC		 	& 248.6 & 53.9 & 479.9 & 318.0 & 1353.5 & 41.4 & 366.8 & 300.5 \\
	\hline
\end{tabular}
\vspace*{0.07cm}
\caption{Results for instances with $n=2000$.
	We report average values over 100 instances each.
  The runtime in seconds is denoted by $t$, while
	$\lambda_1$ denotes the number of vehicles used and
	$\lambda_2$ resp.\
	$\lambda_3$ denote the tour duration resp.\ the travel
        time in hours.}
\label{Tab:2000}
\vspace*{-0.5cm}
\end{table}

In Table \ref{Tab:2000} we present the results of the computational study for $n=2000$.
In the case that the vehicle capacity is the more limiting factor, i.e., $C=200$,
 \nameA\ and \nameC\
 produced nearly identical results with respect to all three objectives.
 In terms of $\lambda_1$ both algorithms clearly  outperform \nameB.
  However, with respect to $\lambda_2$, \nameB\ produced
the best results.
In the case that time window constraints pose the strongest restriction, i.e., $C=400$,
\nameB\ clearly produced  the best results with respect to $\lambda_1$ and $\lambda_2$.
In general, we suggest to apply \nameC\ as an algorithm producing good solutions in most
cases. If the  time window constraints are the most limiting factor of the considered instances,
or if
$\lambda_2$
is more relevant than
$\lambda_1$,
 then \nameB\ is the best choice.

As it is rather hard to remove a whole vehicle from a schedule, we
notice that applying our local improvement heuristic rarely results in a reduction of
the primary objective $\lambda_1$.
However, the study shows that the improvement heuristic heavily impacts
the schedule duration $\lambda_2$, while also having positive impact on the
travel time $\lambda_3$.
The experiments show
that our heuristics,
when applied to large instances containing $2000$
customers,
 take a few seconds to around
five minutes to produce a feasible schedule, and below $25$ minutes to produce an improved
schedule.

\end{document}